\documentclass [12pt,a4paper]{amsart}

\usepackage{amssymb,amsxtra,amsfonts}
\usepackage{lscape}
\newcommand\beq{\begin{equation}}
\newcommand\eeq{\end{equation}}

\newcommand{\IP}{\mathbb{P}}                                     
                           
\newcommand{\IR}{\mathbb{R}}                           
\newcommand{\IC}{\mathbb{C}}
\newcommand{\IZ}{\mathbb{Z}}

\newcommand{\M}{\mathcal{M}}

\newcommand{\cA}{\mathcal{A}}

\newcommand{\cC}{\mathcal{C}}

\newcommand{\bcO}{\boldsymbol{\mathcal{O}}}

\newcommand{\F}{\mathcal{F}}
\newcommand{\cF}{\mathcal{F}}

\newcommand{\cO}{\mathcal{O}}

\newcommand{\g}{       \mathfrak{g}     }

\newcommand{\PVI}{$\text{P}_{\text{VI}}$}   %

\newcommand{\lsl}{\mathfrak{sl}} 

\newcommand{\bA}{{\bf A}}

\newcommand{\pf}{\begin{bpf}}

\newcommand{\pfms}{\begin{bpfms}}
\newcommand{\epf}{\end{bpf}\hfill$\square$\\}           %
\newcommand{\epfms}{\end{bpfms}\hfill$\square$\\}               %

\newcommand{\idea}{\begin{bidea}}

\newcommand{\eidea}{\end{bidea}\hfill$\square$\\}           %

\newcommand{\sk}{\begin{bsk}}    %

\newcommand{\esk}{\end{bsk}\hfill$\square$\\}           %
\newcommand{\sketch}{\begin{bsketch}}%

\newcommand{\esketch}{\end{bsketch}\hfill$\square$\\}

\newcommand{\wt}{\widetilde}

\newcommand{\wh}{\widehat}
\newcommand{\al}{\alpha}
\newcommand{\alc}{\alpha\spcheck}
\newcommand{\be}{\beta}
\newcommand{\ga}{\gamma}
\newcommand{\de}{\delta}

\newcommand{\la}{\lambda}

\newcommand{\om}{\omega}

\newcommand{\res}{{\text{\rm Res}}}

\newcommand{\tr}{\text{\rm Tr}}

\newcommand{\Hom}{\text{\rm Hom}}

\newcommand{\SL}{\text{\rm SL}}

\newcommand{\GL}{\text{\rm GL}}

\newcommand{\End}{\text{\rm End}}
\newcommand{\diag}{{\text{\rm diag}}}

\newcommand{\Aff}{{\text{\rm Aff}}}

\newcommand{\sdp}{{\ltimes}}

\newcommand {\eps}{\varepsilon}

\newcommand{\spq}{/\!\!/}

\def\mapright#1{\smash{
        \mathop{\longrightarrow}\limits^{#1}}}

\theoremstyle{plain}
\newtheorem {hypo}{\bf\hspace{-\parindent}Hypothesis}

\newtheorem {thm}{Theorem}
\newtheorem*{plaincor}{Corollary}

\theoremstyle{definition}

\theoremstyle{remark}
\newtheorem {rmk}[hypo]{Remark}%

\begin{document}

\title[Six results on Painlev\'e VI]
{Six results on Painlev\'e VI}
\author{P. P. Boalch}
\address{\'Ecole Normale Sup\'erieure\\
45 rue d'Ulm\\
75005 Paris\\
France} 
\email{boalch@dma.ens.fr}

\begin{abstract}
After recalling some of
the geometry of the sixth Painlev\'e equation, we will
describe how
the Okamoto symmetries arise naturally from symmetries of Schlesinger's
equations and summarise the classification
of the Platonic Painlev\'e six solutions.
\end{abstract}

\maketitle

\renewcommand{\baselinestretch}{1.02}            %
\normalsize

\begin{section}{Background}

The Painlev\'e VI equation 
is a second order nonlinear differential equation which governs 
the isomonodromic deformations of linear Fuchsian systems of differential 
equations of form
\beq\label{eq: lin syst}
\frac{d}{dz}-\left(
\frac{A_1}{z}+\frac{A_2}{z-t}+\frac{A_3}{z-1}\right),\qquad
A_i\in\g:=\lsl_2(\IC)
\eeq
as the second pole position $t$ varies in $B:=\IP^1\setminus\{0,1,\infty\}$.
(The general case---varying all four pole positions---reduces to this case
using automorphisms of $\IP^1$.)

By `isomonodromic deformation' one means that as $t$ varies 
the linear monodromy representation 
$$\rho:\pi_1(\IP^1\setminus\{0,t,1,\infty\})\to \SL_2(\IC)$$
of \eqref{eq: lin syst} does not change (up to overall conjugation).
Of course this is not quite well-defined since as $t$ varies one is taking
fundamental groups of different four-punctured spheres, and it is crucial to
understand this in order 
to understand the global behaviour (nonlinear monodromy) 
of \PVI\ solutions. 
For small
changes of $t$ there are canonical isomorphisms between the fundamental
groups:
if $t_1,t_2$ are in some disk $\Delta\subset B$ 
in the three-punctured sphere then one has a canonical isomorphism
$$\pi_1(\IP^1\setminus\{0,t_1,1,\infty\})\cong
\pi_1(\IP^1\setminus\{0,t_2,1,\infty\})$$
coming from the homotopy equivalences
$$\IP^1\setminus\{0,t_1,1,\infty\}\hookrightarrow
\{(t,z)\in \Delta\times\IP^1 \bigl\vert z\ne 0,t,1,\infty\}\hookleftarrow
\IP^1\setminus\{0,t_2,1,\infty\}.$$
(Here we view the central space as a family of four-punctured spheres
parameterised by $t\in\Delta$ and are simply saying that it contracts onto any
of its fibres.)

In turn, by taking the space of such $\rho$'s, i.e. the space of
conjugacy classes 
of $\SL_2(\IC)$ representations of the above
fundamental groups one obtains canonical isomorphisms:
$$\Hom(\pi_1(\IP^1\setminus\{0,t_1,1,\infty\}),G)/G\ \cong\ 
\Hom(\pi_1(\IP^1\setminus\{0,t_2,1,\infty\}),G)/G$$
where $G=\SL_2(\IC)$.
Geometrically this says that the spaces of representations
$$\wt M_t:=\Hom(\pi_1(\IP^1\setminus\{0,t,1,\infty\}),G)/G$$
constitute a `local system of varieties' parameterised by $t\in B$. 
In other words, the natural fibration
$$\wt M:=\{(t,\rho)\ \bigl\vert\ t\in B, \rho\in \wt M_t\ \}\mapright{} B$$
over $B$ (whose fibre over $t$ is $\wt M_t$) has a natural flat (Ehresmann)
connection on it.
Moreover this connection is complete:
over any disk in $B$ any two fibres have a canonical identification.  

To get from here to \PVI\ one pulls back the connection on the fibre bundle
$\wt M$ along the Riemann--Hilbert map and writes down the resulting connection
in certain coordinates.
Consequently we see immediately that the monodromy of \PVI\ solutions
corresponds (under the Riemann--Hilbert map)
to the monodromy of the connection on the fibre bundle $\wt M$. 
However since this connection is flat and complete, its monodromy 
is given by the action of the fundamental group of the base $\pi_1(B)\cong
\cF_2$ (the free group on $2$ generators) on a fibre $\wt M_t\subset \wt M$,
which  can easily be written down explicitly, 

Before describing this in more detail 
let us first restrict to linear representations $\rho$ having local
monodromies in fixed conjugacy classes:
$$M_t:=\{ \rho\in \wt M_t\ \bigl\vert\ \rho(\ga_i)\in\cC_i, i=1,2,3,4\}
\subset \wt M_t$$
where $\cC_i\subset G$ are four chosen conjugacy classes, and $\ga_i$ is a
simple positive loop in $\IP^1\setminus\{0,t,1,\infty\}$
around $a_i$, where $(a_1,a_2,a_3,a_4)=(0,t,1,\infty)$ are the four pole
positions. 
(By convention we assume 
the loop $\ga_4\cdots\ga_1$ is contractible, and 
note that $M_t$ is two-dimensional in general.)
The connection on $\wt M$ restricts to a (complete flat  Ehresmann) 
connection on the fibration
$$M:=\{(t,\rho)\ \bigl\vert\ t\in B, \rho\in M_t\ \}\to B$$
whose fibre over $t\in B$ is $M_t$.
The action of $\F_2=\pi_1(B)$ on the fibre $M_t$ 
(giving the monodromy of the connection on the bundle $M$
and thus the monodromy of the corresponding \PVI\ solution)
is given explicitly as follows.
Let $w_1,w_2$ denote the generators of $\F_2$, thought of as simple positive
loops in $B$ based at $1/2$ encircling $0$ (resp. $1$) once.
Then $w_i$ acts on $\rho\in M_t$ as the square of $\om_i$ where
$\om_i$ acts by fixing $M_j$ for $j\ne i,i+1, \ (1\leqslant j\leqslant 4)$
and
\beq \label{eq: f2 action}
\om_i(M_i,\ M_{i+1})= (M_{i+1},\ M_{i+1}M_{i}M^{-1}_{i+1})
\eeq
where $M_j=\rho(\ga_j)\in G$ is the $j$th monodromy matrix.
Indeed $\F_2$ can naturally be identified with the pure mapping class group of
the four-punctured sphere and this action comes from its natural action
(by push-forward of loops) as
outer automorphisms of $\pi_1(\IP^1\setminus\{0,t,1,\infty\})$, 
cf. \cite{Birman}.
(The geometric origins of this action in the context of isomonodromy can be
traced back at least to Malgrange's work \cite{Mal-imd1long} 
on the global properties of the Schlesinger equations.)

On the other side of the Riemann--Hilbert correspondence we may choose some
adjoint orbits $\cO_i\subset\g:=\lsl_2(\IC)$
such that $$\exp(2\pi\sqrt{-1}\cO_i)=\cC_i$$
and construct the space of residues:
$$\bcO:= \cO_1\times\cdots \times\cO_4\spq G=
\left\{\ (A_1,\ldots A_4)\in \cO_1\times \cdots\times\cO_4 
\ \bigl\vert \ \text{$\sum A_i =0$}\ \right\}/G$$ 
where, on the right-hand side, $G$ is acting by diagonal conjugation:
$g\cdot(A_1,\ldots A_4)=(gA_1g^{-1},\ldots,gA_4g^{-1})$.
This space $\bcO$ is also two-dimensional in general.
To construct a Fuchsian system \eqref{eq: lin syst}
out of such a four-tuple of
residues one must also choose a value of $t$, so the total space of
linear connections we are interested in is:
$$\M^*:=\bcO\times B$$
and we think of a point $(\bA,t)\in \M^*$, where $\bA=(A_1,\ldots,A_4)$,
as representing the linear connection
$$\nabla=d-Adz,\qquad \text{where }A=\sum_1^3\frac{A_i}{z-a_i},
\quad(a_1,a_2,a_3,a_4)=(0,t,1,\infty)$$
or equivalently the Fuchsian system \eqref{eq: lin syst}.
(Observe that $A_i=\res_{a_i}(Adz)$.) 

If we think of $\M^*$ as being a (trivial) fibre bundle over $B$ with fibre 
$\bcO$ then, provided the residues are sufficiently generic (e.g. if no
eigenvalues differ by positive integers), the Riemann--Hilbert map (taking
linear connections to their monodromy representations) gives a bundle map
$$\nu:\M^*\to M.$$ 
Written like this the Riemann--Hilbert map $\nu$ is a holomorphic
map (which is in fact injective if the eigenvalues are also nonzero
cf. e.g. \cite{JMU81} Proposition 2.5).
We may then pull-back (restrict) the nonlinear connection on $M$ to give a
nonlinear connection on the bundle $\M^*$, which we will refer to as the 
{\em isomonodromy connection}.

The remarkable fact is that even though the Riemann--Hilbert map is
transcendental, the connection one obtains in this way is algebraic.
Indeed Schlesinger \cite{Schles} showed that locally 
horizontal sections $\bA(t):B\to \M^*$ are given
(up to overall conjugation) by solutions to   
the Schlesinger equations:
\begin{equation} \label{eqn: 2x2 schles}
\frac{dA_1}{dt}=\frac{[A_2,A_1]}{t},\qquad 
\frac{dA_2}{dt}=\frac{[A_1,A_2]}{t}+\frac{[A_3,A_2]}{t-1},
\qquad 
\frac{dA_3}{dt}=\frac{[A_2,A_3]}{t-1} 
\end{equation}
which are (nonlinear) {\em algebraic} differential equations. 

To get from the Schlesinger equations to \PVI\ one proceeds as follows 
(cf. \cite{JM81} Appendix C).
Label the eigenvalues of $A_i$ by $\pm\theta_i/2$ (thus choosing an order of
the eigenvalues or equivalently, if the reader prefers, a quasi-parabolic
structure at each singularity), and suppose $A_4$ is diagonalisable.
Conjugate the system so that 
$$A_4=-(A_1+A_2+A_3)=\diag(\theta_4,-\theta_4)/2$$
and note that Schlesinger's equations preserve $A_4$.
Since the top-right matrix entry of $A_4$ is zero, the top-right
 matrix entry of 
\beq\label{eq: matrix poly}
z(z-1)(z-t)\sum_1^3\frac{A_i}{z-a_i}
\eeq
is a degree one polynomial in $z$. Define $y(t)$ to be the position
of its unique zero on the complex $z$ line.

\begin{thm}[see \cite{JM81}]\label{thm: 2x2 schles to P6}
If $\bA(t)$ satisfies the Schlesinger equations then $y(t)$ satisfies \PVI:
\begin{align*}\frac{d^2y}{dt^2}=
&\frac{1}{2}\left(\frac{1}{y}+\frac{1}{y-1}+\frac{1}{y-t}\right)
\left(\frac{dy}{dt}\right)^2
-\left(\frac{1}{t}+\frac{1}{t-1}+\frac{1}{y-t}\right)\frac{dy}{dt}\\
&+\frac{y(y-1)(y-t)}{2\,t^2(t-1)^2}\left(
(\theta_4-1)^2-
\frac{\theta_1^2\, t}{y^2}+ 
\frac{\theta_3^2(t-1)}{(y-1)^2}+
\frac{(1-\theta_2^2)t(t-1)}{(y-t)^2}\right).
\end{align*}
\end{thm}

Phrased differently, for each fixed $t$, the prescription above defines a
function $y$ on $\bcO$, which makes up half of a system of (canonical)
coordinates, defined on a dense open subset.
A conjugate coordinate $x$ can be explicitly 
defined and
one can write the isomonodromy connection explicitly in the coordinates
$x,y$ on $\bcO$ to obtain a coupled system of first-order nonlinear equations
for $x(t), y(t)$ 
(see \cite{JM81}, where our $x$ is denoted $\wt z$). 
Then eliminating $x$ yields the second order equation
\PVI\ for $y$.
(One consequence is that if $y$ solves \PVI\ there is a direct relation
between $x$ and the derivative $y'$, as in equation \eqref{eq: p eq} below.)

In the remainder of this article the main aims are to:

$\bullet$1) Explain how Okamoto's affine $F_4$ Weyl group 
symmetries of \PVI\ arise from
natural symmetries of Schlesinger equations, and

$\bullet$2) Describe the classification of the Platonic solutions 
to \PVI\ (i.e. those
solutions having linear monodromy group equal to the symmetry group
of a Platonic solid).
\vspace{0.2cm}

The key step for $\bullet$1) (which also led us to $\bullet$2))
is to use a different realisation of 
\PVI, as controlling isomonodromic deformations of certain 
$3\times 3$ Fuchsian systems.  
Note that these results have been written down elsewhere,
although 
the explicit formulae of Remarks \ref{rmk: new pf} and
\ref{rmk: simpler Bis} are new 
and constitute a direct verification of 
the main results about the $3\times 3$ realisation.
Note also that 
the construction of the Platonic solutions has 
evolved rapidly recently (e.g. since the author's talk in Angers
and since the 
first version of \cite{icosa}
appeared). For example there are now simple explicit formulae for 
all the Platonic
solutions (something that we had not imagined was 
possible for a long time\footnote{
Mainly because the 18 branch genus one icosahedral solution of \cite{DubMaz00}
took 10 pages to write down and we knew quite early on that the
largest icosahedral solution had genus seven and 72 branches.}).

\begin{rmk}
Let us briefly mention some other possible directions that will
not be discussed further here.
Firstly, by 
describing \PVI\  
in this way the author is trying to emphasise that \PVI\  is the
explicit form of the simplest non-abelian Gauss--Manin connection, in the sense
of Simpson \cite{Sim94ab}, thereby putting \PVI\ 
in a very general context (propounded further in \cite{smid} section 7, 
especially p.192). 
For example suppose we replace the above
family of four-punctured spheres (over $B$)
by a family of projective varieties $X$ over a base $S$, and choose a
complex reductive group $G$.
Then (by the same argument as above) one again has a local system of varieties
$$M_B= \Hom(\pi_1(X_s),G)/G$$ over $S$ 
and one can pull-back along the Riemann--Hilbert map to obtain a flat
connection on 
the corresponding family $M_{DR}$ of moduli spaces of connections. 
Simpson proves this connection is again algebraic, and calls it the
non-abelian Gauss--Manin connection, since $M_B$ and $M_{DR}$ are two
realisations of the first non-abelian cohomology group $H^1(X_s,G)$, the Betti
and De Rham realisations.

Also, much 
of the structure found in the regular (-singular) case may be generalised
to the irregular case.
For example as Jimbo--Miwa--Ueno \cite{JMU81} showed, one can
also consider isomonodromic deformations of (generic) irregular connections 
on a Riemann surface and obtain explicit deformation equations in the case of
$\IP^1$. This can also be  described in terms of nonlinear 
connections on moduli spaces and there are natural symplectic
structures on the moduli 
spaces which are preserved by the connections 
\cite{smid,saqh}. Perhaps most interestingly one obtains extra deformation
parameters in the irregular case (one may vary the `irregular type' of the
linear connections as well as the moduli of the punctured curve).
These extra deformation parameters 
turn out to be related to quantum Weyl groups \cite{bafi}.

As another example, in the regular (-singular) case
non-abelian Hodge theory \cite{Simpson-nabh}
gives a third 
``Dolbeault'' realisation of $H^1(X_s,G)$ as a moduli space of Higgs bundles, 
closely related to the 
existence of a  hyperK\"ahler structure on the moduli space.
The moduli spaces of (generic) irregular connections
on curves may also be realised in terms of Higgs bundles and admit
hyperK\"ahler metrics \cite{nabh}.
\end{rmk}

\end{section}

\begin{section}{Affine Weyl group symmetries}
If we subtract off $y''=\frac{d^2y}{dt^2}$ from the right-hand side of 
the \PVI\ equation and multiply through by 
$t^2(t-1)^2y(y-1)(y-t)$ then we obtain a polynomial:
$$P(t,y,y',y'', \theta)\in 
\IC[t,y,y',y'',\theta_1,\theta_2,\theta_3,\theta_4]$$
where $\theta=(\theta_1,\theta_2,\theta_3,\theta_4)$ are the parameters.

Suppose $\Pi$ is a Riemann surface equipped with a holomorphic map 
$t:\Pi\to U$ onto some open subset $U\subset B:=\IP^1\setminus\{0,1,\infty\}$,
with non-zero derivative (so $t$ is always a local isomorphism).
(For example one could take $\Pi=U$ with $t$ the inclusion,
or take $\Pi$ to be the upper half-plane, and $t$ the universal covering map
onto $U=B$.)
Then a meromorphic function $y$ on $\Pi$ will be said to be a solution to
\PVI\ if
\beq\label{eq: Pzero}
P(t,y, y',y'', \theta)=0
\eeq
as functions on $\Pi$, 
for some choice of $\theta$, where $y'=\frac{dy}{dt},y''=\frac{d^2y}{dt^2}$ are
defined by using
$t$ as a local parameter on $\Pi$. (With this $t$-dependence understood we will
abbreviate \eqref{eq: Pzero} as $P(t,y, \theta)=0$ below.)
By definition the finite branching solutions to \PVI\ are those with
$\Pi$ a finite cover of $B$, i.e. so that $t$ is a Belyi map.
Such $\Pi$ admits a natural compactification $\overline\Pi$, 
on which $t$ extends to a
rational function.
The solution is 
``algebraic'' if $y$ is a rational function on $\overline \Pi$.
Given an algebraic solution ($\overline\Pi,y,t$) we will say the curve 
$\overline \Pi$ is ``minimal'' or is an ``efficient parameterisation'' if $y$
generates the function field of $\overline \Pi$ as an extension of $\IC(t)$.
The ``degree'' (or number of ``branches'') of an algebraic  solution 
is the degree of the map $t:\overline \Pi\to \IP^1$ (for 
$\overline \Pi$ minimal) and
the genus of the solution is the genus of the (minimal) curve $\overline\Pi$.
(The genus can easily be computed in terms of the
nonlinear monodromy of the \PVI\ solution using the Riemann--Hurwitz formula,
i.e. in terms of the explicit $\cF_2$ action above 
on the linear monodromy data.)

Four symmetries of \PVI\ 
(which we will label $R_1,\ldots,R_4$) are immediate:
\begin{align}
\tag{$R_1$}P(t,y, \theta)&=P(t,y,-\theta_1,\theta_2,\theta_3,\theta_4)\\
\tag{$R_2$}&=P(t,y,\theta_1,-\theta_2,\theta_3,\theta_4)\\
\tag{$R_3$}&=P(t,y,\theta_1,\theta_2,-\theta_3,\theta_4)\\
\tag{$R_4$}&=P(t,y,\theta_1,\theta_2,\theta_3,2-\theta_4)
\end{align}
since $P$ only depends on the squares of 
$\theta_1,\theta_2,\theta_3$ and $\theta_4-1$.

Okamoto \cite{OkaPVI} proved there are also much less trivial symmetries:
\begin{thm}\label{thm: ok.sym}
If $P(t,y, \theta)=0$ then
\begin{align}
\tag{$R_5$}P(t,y+\delta/x,
\theta_1-\delta,\theta_2-\delta,\theta_3-\delta,\theta_4-\delta)&=0
\end{align}
where $\delta=\sum_1^4\theta_i/2$ and
\beq\label{eq: p eq}
2 x=
\frac {\left( t-1 \right)y'-\theta_1}{y}+
\frac{y'-1-\theta_2}{y-t}-
\frac {t\,y'+\theta_3}{y-1}.
\eeq
\end{thm}

\begin{rmk}
This can be verified directly 
by a symbolic computation in differential algebra.
On actual solutions however 
it is not always well-defined since for example one may have
$y=t$ (identically) or find $x$ is identically zero. It seems one can avoid
these problems by assuming $y$ is not a Riccati solution 
(cf. \cite{watanabePVI}). For example if one finds $x=0$ then we see $y$ 
solves a first order (Riccati) equation, so was a Riccati solution.
Moreover the Riccati solutions are well understood and correspond to the
linear representations $\rho$ which are either reducible or rigid, so little
generality is lost.
\end{rmk}

\begin{rmk}
In terms of the symmetries $s_0,\ldots,s_4$ of \cite{NY-so8},
$R_1,\ldots R_4$ are $s_4,s_0,s_3,s_1$ respectively and $R_5$ is 
conjugate to $s_2$ via $R_1R_2R_3R_4$, where the parameters
$\al_4,\al_0,\al_3,\al_1$
of \cite{NY-so8} are taken to be $\theta_1,\theta_2,\theta_3,\theta_4-1$
respectively, and 
$p=x+\sum_1^3\theta_i/(y-a_i)$.  
\end{rmk}

A basic observation (of Okamoto) 
is that these five symmetries generate a group
isomorphic to the affine Weyl group of type $D_4$. 
More precisely let $\eps_1,\ldots,\eps_4$ be an orthonormal basis of a
Euclidean vector space $V_\IR$ with inner product $(\ ,\,)$ 
and complexification 
$V$, and consider the
following set of $24$ unit vectors 
$$D_4^-=\{\pm\eps_i,\  (\pm\eps_1\pm\eps_2\pm\eps_3\pm\eps_4)/2\}.$$
This is a root system isomorphic to the standard $D_4$ root system
$$D_4=\{\pm\eps_i\pm\eps_j (i<j)\}$$
but with vectors of length $1$ rather than $\sqrt{2}$.
(Our main reference for root systems etc. is \cite{BbkLie}.
One may identify $D_4^-$ with the {\em group} of units of the Hurwitzian
integral quaternions \cite{conway-smith}, and then identify with $D_4$ by
multiplying by the quaternion $1+{\bf i}$.)
Each root $\alpha\in D^-_4$ determines a coroot 
$\alc = \frac{2\al}{(\al,\al)}$($=2\al$ here) as well as
a hyperplane $L_\al$ in $V$:
$$L_\alpha:= \{\  v\in V\ \bigl\vert\ (\alpha,v)=0\ \}.$$

In turn $\al$ 
determines an orthogonal reflection $s_\alpha$, the reflection in this
hyperplane:
$$s_\alpha(v)= v- 2\frac{(\al,v)}{(\al,\al)}\al=v- (\alc,v)\al.$$

The Weyl group $W(D^-_4)\subset O(V)$ is the group 
generated by these reflections:
$$W(D^-_4)= \langle\ s_\al\ \bigl\vert\ \al\in D^-_4\ \rangle$$
which is of order $192$.
Similarly the choice of a root $\alpha\in D^-_4$ and an integer $k\in \IZ$
determines an affine hyperplane $L_{\al,k}$ in $V$:
$$L_{\alpha,k}:= \{\  v\in V\ \bigl\vert\ (\alpha,v)=k\ \}$$
and the reflection $s_{\alpha,k}$ 
in this hyperplane is an affine Euclidean transformation 
$$s_{\al,k}(v)= s_\al(v)+k\alc.$$

The affine Weyl group $W_a(D^-_4)\subset \Aff(V)$ is the group 
generated by these reflections:
$$W_a(D^-_4)= 
\langle\ s_{\al,k}\ \bigl\vert\ \al\in D^-_4,\  k\in\IZ\ \rangle$$
which is an infinite group isomorphic to the semi-direct product
of $W(D^-_4)$ and the coroot lattice $Q((D^-_4)\spcheck)$ 
(which is the lattice in
$V$ generated by the coroots $\alc\in(D^-_4)\spcheck=D_4^+=2D_4^-$).
By definition the connected components of the complement in $V_\IR$ 
of all the (affine) reflection
hyperplanes are the $D^-_4$ alcoves. 
The closure $\overline \cA$ 
in $V_\IR$ of any alcove $\cA$ is a fundamental domain for the
action of the affine Weyl group; every $W_a(D_4^-)$ orbit in $V_\IR$ 
intersects
$\overline \cA$ in precisely one point.

Now if we write a point of $V$ as $\sum \theta_i\eps_i$ (i.e. the parameters
$\theta_i$ are being viewed as coordinates on $V$ with respect to the 
$\eps$-basis) then, on $V$, 
the five symmetries above correspond to the reflections in the 
five hyperplanes:
$$\theta_1=0,\qquad\theta_2=0,\qquad\theta_3=0,\qquad\theta_4=1,
\qquad\sum\theta_i=0.$$
The reflections in these hyperplanes 
generate $W_a(D_4^-)$ since the region:
$$\theta_1<0,\qquad\theta_2<0,\qquad\theta_3<0,
\qquad\theta_4<1,\qquad\sum\theta_i>0$$
that they bound in $V_\IR$ is an alcove. 
(With respect to the root ordering 
given by taking the inner product with the vector $4\eps_4-\sum_1^3\eps_i$, 
the roots
$-\eps_1,-\eps_2,-\eps_3,\sum\eps_i/2$ are a basis of positive roots of
$D_4^-$, and the highest root is $\eps_4$, so by \cite{BbkLie} (p.175) 
this is an alcove.)

In fact, as Okamoto showed, the full symmetry group of \PVI\ is the affine
Weyl group of type $F_4$. 
(The $F_4$ root system is the set of  $48$ vectors in the
union of $D_4$ and $D_4^-$.)
This is not surprising if one recalls that
$W_a(F_4)$ is the normaliser of $W_a(D_4^-)$ in the group of affine
transformations; $W_a(F_4)$ is the extension of $W_a(D_4^-)$ by
the symmetric group on four letters, $S_4$ thought of 
as the automorphisms of the affine 
$D_4$ Dynkin diagram (a central node with four satellites).
This extension breaks into two pieces corresponding to the exact sequence
$$1\longrightarrow K_4 \longrightarrow S_4 
\longrightarrow S_3 \longrightarrow 1$$
where $K_4\cong (\IZ/2)^2$  is the Klein four-group.
On one hand the group of translations 
is extended by a $K_4$; the lattice
$Q(D_4^+)$ is replaced by $Q(F_4\spcheck)=Q(D_4)$.
(In general \cite{BbkLie} p.176 
one replaces $Q(R\spcheck)$ by $P(R\spcheck)=Q(R)^*$.)
On the other hand the Weyl group is extended by an $S_3$, thought of as the
automorphisms of the usual $D_4$ Dynkin diagram;
$W(D^-_4)$ is replaced by the full group of automorphisms $A(D^-_4)$ of the 
root system, which in this case is equal to $W(F_4)$.

Likewise the corresponding symmetries of \PVI\ break into two pieces.
First one has an $S_3$ permuting $\theta_i$ ($i=1,2,3$) generated for example
by the symmetries (denoted $x^1,x^3$ respectively in \cite{OkaPVI} p.361):

\begin{align*}
P(t,y, \theta)=0 \quad \Longrightarrow\quad
&P(1-t,1-y,\theta_3,\theta_2,\theta_1,\theta_4)=0\\
P(t,y, \theta)=0\quad \Longrightarrow\quad
&P\left(\frac{t}{t-1},\frac{t-y}{t-1},\theta_2,\theta_1,\theta_3,\theta_4
\right)=0.
\end{align*}

We remark that $W_a(D_4^-)$ already contains transformations permuting $\theta$
by the standard Klein four group (mapping $\theta$ to 
$(\theta_3,\theta_4,\theta_1,\theta_2)$ etc.), and so we already obtain all
permutations of $\theta$ just by adding the above two symmetries.\footnote{
For example $R_5r_1r_3R_5r_2r_4$ produces the permutation written, where
$r_i$ is the Okamoto transformation negating $\theta_i$---i.e.
$r_i=R_i$ for $i=1,2,3$ and $r_4=R_5(R_1R_2R_3)R_5(R_1R_2R_3)R_5$.} 
To obtain the desired $K_4$ extension we refine the possible translations by
adding the further symmetry (denoted $x^2$ in \cite{OkaPVI}):
$$
P(t,y, \theta)=0\quad \Longrightarrow\quad
P(1/t,1/y,\theta_4-1,\theta_2,\theta_3,\theta_1+1)=0.
$$
Combined with $x^1,x^3$ this generates an $S_4$ which may be thought of as
permuting the set of values of $\theta_1,\theta_2,\theta_3,\theta_4-1$.
(Note that, modulo the permutations of $\theta$, we now have 
translations of the form
$\theta\mapsto (\theta_1+1,\theta_2,\theta_3,\theta_4-1)$, generating 
$Q(D_4)$.)
\begin{rmk}
One can also just extend by the $K_4$ and get an intermediate group, often
called the extended Weyl group $W'_a(D_4^-)=W(D_4^-)\sdp P((D_4^-)\spcheck)$
which is normal in $W_a(F_4)$ and is the maximal subgroup that does not change 
the time $t$ in the above action on \PVI. The quotient group $S_3$ should thus
be thought of as the automorphisms of $\IP^1\setminus\{0,1,\infty\}$.
\end{rmk}

Our aim in the rest of this section is to explain how these symmetries arise
naturally from symmetries of the Schlesinger equations.
The immediate symmetries are: 

\ 

$\bullet$\ 
(twisted) Schlesinger transformations,

$\bullet$\ negating the $\theta_i$ independently, and 

$\bullet$\ arbitrary permutations of the $\theta_i$.

\ 

In more detail the Schlesinger transformations (see \cite{JM81}) are certain
rational gauge transformations which shift the eigenvalues of the residues by
integers. Applying such a transformation and then twisting by a logarithmic
connection on the trivial line bundle (to return the system to $\lsl_2$)
is a symmetry of the Schlesinger equations.
(This procedure of ``twisting'' clearly commutes with the flows of the
Schlesinger equations: in concrete terms it 
simply amounts to adding an expression 
of the form $\sum_1^3 c_i/(z-a_i)$, for constant scalars $c_i$, 
to the Fuchsian system
\eqref{eq: lin syst}. Recall $(a_1,a_2,a_3)=(0,t,1)$.)

Secondly the eigenvalues of the residues are only determined by the abstract 
Fuchsian 
system up to sign (i.e. one chooses an order of the eigenvalues of each
residue to define $\theta_i$, and these choices can be swapped).

Finally if we permute the labels $a_1,\ldots,a_4$ of 
the singularities of the Fuchsian 
system arbitrarily and then perform the (unique) 
automorphism of the sphere mapping
$a_1,a_3,a_4$ to $0,1,\infty$ respectively, we obtain another isomonodromic
family of systems, which can be conjugated to give another 
Schlesinger solution.

As an example consider the case of negating $\theta_4$.
Suppose we have a solution of the Schlesinger equations 
$\bA(t)$ for a given choice of $\theta$
and have normalised $A_4$ as required in 
Theorem \ref{thm: 2x2 schles to P6} (this is where the sign choice is
used). If we conjugate $\bA$ by the permutation matrix
$\left(\begin{smallmatrix}0&1\\1&0\end{smallmatrix}\right)$
we again get a solution of the Schlesinger equations, and by Theorem
\ref{thm: 2x2 schles to P6} this yields a solution to \PVI\ with parameters
$(\theta_1,\theta_2,\theta_3,-\theta_4)$.
This gives the corresponding Okamoto transformation in terms of Schlesinger
symmetries. (It is a good, if unenlightening, exercise 
to compute the explicit 
formula---in effect computing the position of the zero of the bottom-left entry
of \eqref{eq: matrix poly} in terms of $x,y$---and check it agrees with the
action of the 
corresponding word in the given generators of $W_a(D_4^-)$, although logically
this verification is unnecessary since a) This is a symmetry of \PVI\ and b)
Okamoto found all symmetries, and they are determined by their action on
$\{\theta\}$.)

However one easily sees
that the group generated by these immediate symmetries 
does not contain the
transformation $R_5$ of Theorem \ref{thm: ok.sym}.
To obtain this symmetry
we will recall (from \cite{k2p})  
how \PVI\ also governs the 
isomonodromic deformations of certain rank three Fuchsian systems 
and show that $R_5$
arises from symmetries of the corresponding Schlesinger equations 
(indeed it arises simply from the choice of ordering of the 
eigenvalues at infinity).
(Note that Noumi--Yamada \cite{NY-so8} have also obtained this symmetry from an
isomonodromy viewpoint, but only in terms of
an {\em irregular} (non-Fuchsian) $8\times 8$ system whose
isomonodromy deformations, in a generalised sense, are governed by \PVI.)

To this end
let $V=\IC^3$ be a three-dimensional complex vector space
and suppose $B_1,B_2,B_3\in \End(V)$ are rank one matrices.
Let $\lambda_i=\tr(B_i)$ and suppose that $B_1+B_2+B_3$ is diagonalisable with 
eigenvalues $\mu_1,\mu_2,\mu_3$, so that taking the trace implies
\beq\label{eq: suml eq summ}
\sum_1^3\lambda_i=\sum_1^3\mu_i.
\eeq
Consider connections of the form
\beq\label{eq: 3x3 conn}
\nabla=d-\wh Bdz, \qquad 
\wh B(z)=\frac{B_1}{z}+\frac{B_2}{z-t}+\frac{B_3}{z-1}.
\eeq

The fact is that the isomonodromic deformations of such connections are also
governed by \PVI\ (one might expect such a thing since the corresponding
moduli spaces are again two-dimensional). 
One proof of this (\cite{pecr}) 
is to show directly that the corresponding Schlesinger
equations are equivalent to those arising in the original $2\times 2$ 
case (this may done easily by writing out the isomonodromy connections 
explicitly in
terms of the coordinates on the spaces of residues 
given by the invariant functions, and
comparing the resulting nonlinear differential equations).

The second proof of this result directly gives the function that
solves \PVI;
First conjugate $B_1,B_2,B_3$ by a single element of $\GL_3(\IC)$ such that 
$$B_1+B_2+B_3=\diag(\mu_1,\mu_2,\mu_3).$$
(Note this uses the choice of ordering of eigenvalues of $B_1+B_2+B_3$.)
Consider the polynomial defined to be the $(2,3)$ matrix entry of
\beq\label{eq: 3x3 poly}
z(z-1)(z-t)\wh B(z).\eeq
By construction this is a linear polynomial, so has a unique zero on the
complex plane. Define $y=y_{23}$ to be the position of this zero.

\begin{thm}[\cite{k2p} p.201]\label{thm: direct PVI}
If we vary $t$ and evolve $\wh B$ according to Schlesinger's equations 
then $y(t)$ satisfies the \PVI\  equation 
with parameters 
\begin{equation}\label{thla}
\theta_1= \lambda_1-\mu_1,\quad
\theta_2= \lambda_2-\mu_1,\quad
\theta_3= \lambda_3-\mu_1,\quad
\theta_4=\mu_3-\mu_2.
\end{equation}
\end{thm}

The proof given in \cite{k2p} uses an extra symmetry of the corresponding
Schlesinger equations (\cite{k2p} Proposition 16) to pass to the $2\times 2$
case. Note that \cite{k2p} also gives the explicit 
relation between the $2\times 2$ and $3\times 3$ linear monodromy data, not
just the relation between the Fuchsian systems. 

\begin{rmk}
Apparently (\cite{Dett-Reit-p6}) 
this procedure of \cite{k2p} is essentially
N. Katz's middle-convolution functor \cite{Katz-rls} in this context. 
For us it originated by considering the effect of performing the
Fourier--Laplace transformation, twisting by a flat line bundle $\lambda dw/w$
and transforming back (reading \cite{BJL79} carefully 
to see what happens to the connections and their monodromy). 
It is amusing that the middle-convolution functor first arose
through the $l$-adic Fourier transform, 
essentially in this way it seems, 
and was then translated back into the complex analytic world, rather than
having been previously worked out directly.
\end{rmk}

If we now conjugate $\wh B(z)$ by an arbitrary $3\times 3$ permutation matrix
(i.e. a matrix which is zero except for 
precisely one $1$ in each row and column),
we obtain another solution of the Schlesinger equations, but with the
$\mu_i$ permuted accordingly. The happy fact that this 
$S_3$ transitively permutes the six off-diagonal entries yields:
\begin{plaincor}
Let $(i,j,k)$ be some permutation of $(1,2,3)$. Then the position $y_{jk}$
of the zero of the $(j,k)$ matrix entry of \eqref{eq: 3x3 poly} satisfies
\PVI\ with parameters
\begin{equation}\label{thlaijk}
\theta_1= \lambda_1-\mu_i,\quad
\theta_2= \lambda_2-\mu_i,\quad
\theta_3= \lambda_3-\mu_i,\quad
\theta_4=\mu_k-\mu_j.
\end{equation}
\end{plaincor}
\pf Conjugate by the corresponding permutation matrix and apply Theorem
\ref{thm: direct PVI}. \epf 

For example the permutation swapping $\mu_2$ and $\mu_3$
 thus amounts to negating $\theta_4$ (indeed
one may view the original $2\times 2$ picture as embedded in this $3\times 3$
picture as the bottom-right $2\times 2$ 
submatrices, at least after 
twisting by a logarithmic connection on a line bundle to make 
$A_1,A_2,A_3$ rank one matrices).

More interestingly let us compute the action on the
$\theta$ parameters of the permutation swapping
$\mu_1$ and $\mu_3$:
\begin{align*}
\theta&= (
\lambda_1-\mu_1,
\lambda_2-\mu_1,
\lambda_3-\mu_1,
\mu_3-\mu_2),\\
\theta'&= (
\lambda_1-\mu_3,
\lambda_2-\mu_3,
\lambda_3-\mu_3,
\mu_1-\mu_2).
\end{align*}
Thus $\theta'_i=\theta_i-\delta$ with $\delta=\mu_3-\mu_1$.
However using the relation \eqref{eq: suml eq summ}
we find 
$$\sum_1^4 \theta_i=
\sum_1^3\lambda_i-3\mu_1+\mu_3-\mu_2=
2(\mu_3-\mu_1)$$
so that $\delta=\sum_1^4 \theta_i/2$ as required for $R_5$. This leads to:
\begin{thm}[\cite{k2p} p.202]\label{thm: R5 from schles}
The permutation swapping $\mu_1$ and $\mu_3$
yields the Okamoto transformation $R_5$.
In other words if $y=y_{23}$ and $\delta=\sum_1^4 \theta_i/2$ and 
$$2 x=
\frac {\left( t-1 \right)y'-\theta_1}{y}+
\frac{y'-1-\theta_2}{y-t}-
\frac {t\,y'+\theta_3}{y-1}$$
then
$$y_{21}=y+\frac{\delta}{x}.$$
\end{thm}

\begin{rmk} \label{rmk: new pf}
Of course if one had a suitable parameterisation of the space of such 
$3\times3$  linear
connections \eqref{eq: 3x3 conn}
in terms of $x$ and $y$, this could be proved by a direct computation.
Such a parameterisation may be obtained as follows 
(lifted from the $2\times 2$ case in \cite{JM81} using \cite{k2p} Prop. 16).
(In particular this shows how one might have obtained the transformation 
formula of Theorem \ref{thm: ok.sym} directly.)
Fix $\la_i,\mu_i$ for $i=1,2,3$ such that $\sum\la_i=\sum\mu_i$.
We wish to write down the matrix entries of $B_1,B_2,B_3$ as rational functions
of $x,y,t,\la_i,\mu_i$. 
The usual $2\times 2$ parameterisation of Jimbo--Miwa \cite{JM81}
will appear in the bottom-right
corner if $\mu_1=0$.
First define $\theta_i$ as in Theorem \ref{thm: direct PVI}.
Then define $z_i,u_i$ for 
$i=1,2,3$ as the unique solution to the $6$ equations:
$$y=tu_1z_1,\qquad x=\sum z_i/(y-a_i),\qquad\sum z_i=\mu_1-\mu_3,$$
$$\sum u_iz_i=0,\qquad\sum w_i=0,\qquad\sum(t-a_i)u_iz_i=1,$$
where $w_i=(z_i+\theta_i)/u_i$ and $(a_1,a_2,a_3)=(0,t,1)$
(cf. \cite{JM81} and \cite{octa} Appendix A).
Now define $c_1,c_2,c_3$ as the solution to the $3$ linear equations:
$$\sum c_iz_i=0,\quad\sum c_iw_i=0,\quad\sum(t-a_i)c_iz_i=1.$$
The determinant of the corresponding $3\times 3$ 
matrix is generically nonzero so this
yields explicit formulae for the $c_i$ (using for example the formula
for the inverse of a $3\times 3$ matrix)---we will not write them since they
are somewhat clumsy and easily derived from the above equations.\footnote{
For the reader's convenience a text file with some Maple code to verify 
the
assertions of this remark (and some others in this article) is available 
at www.dma.ens.fr/\~\ \!\!boalch/files/sps.mpl (or alternatively 
with the source file of arxiv:math.AG/0503043).}
Using $z_i,u_i,w_i,c_i$  we construct  
forms $\be_i$ and vectors $f_i$ for $i=1,2,3$  by setting
$$\be_i=(0 , w_i, -z_i)\in V^*,\qquad 
f_i=\left(\begin{matrix}c_i\\u_i\\1 \end{matrix}\right)\in V.$$
(The meaning of the above $9$ equations is simply that 
if we set $B_i^0=f_i\otimes \be_i\in \End(V)$ and
$\wt B^0=z(z-1)(z-t)\wh B^0$ where $\wh B^0=\sum B_i^0/(z-a_i)$
then 
$$\sum B_i^0=\diag(\mu_1,\mu_2,\mu_3)-\mu_1,\qquad
-\wh B^0_{33}\bigl\vert_{z=y}=x,\qquad
\wt B^0_{23}=z-y$$
and the coefficient of $z$ in  the top-right entry 
$\wt B^0_{13}$ is also $1$.)

The $f_i$ are in general linearly independent and we can define
the dual basis $\wh f_i\in V^*$, with $\wh f_i(f_j)=\delta_{ij},$ explicitly.
The desired matrices are then 
$$B_i=f_i\otimes(\be_i+\mu_1\wh f_i)\in \End(V).$$ 
Clearly $B_i$ is a rank-one matrix and one may check that $\tr(B_i)=\la_i$
and that $\sum B_i=\diag(\mu_1,\mu_2,\mu_3)$.
Moreover generically any such triple 
of rank-one matrices is conjugate to the triple
$B_1,B_2,B_3$
up to overall conjugation by the diagonal torus, for some values of $x$ and 
$y$.
Now if we define $y_{ij}$ to be the value of $z$ for which
the $i,j$ matrix entry of $\wt B:=z(z-1)(z-t)\wh B$ vanishes, where
$\wh B=\sum_1^3B_i/(z-a_i)$
then one may check explicitly (e.g. using Maple) that $y_{23}=y$
and $y_{21}=y+(\mu_3-\mu_1)/x$ as required.  
Also $x$ may be defined in general, 
as a function on the space of such connections, by the
prescription:
$$x= \frac{\mu_1-\mu_3}{\mu_3}\wh B_{33}\bigl\vert_{z=y}$$
which may be checked to hold in the above parameterisation, and specialises to
the usual definition of $x$ in the $2\times 2$ case when $\mu_1=0$.
Moreover, one may check $x$ is preserved under $R_5$, and this agrees with the
fact that one also has
$$x= \frac{\mu_3-\mu_1}{\mu_1}\wh B_{11}\bigl\vert_{z=y+\delta/x}$$
in the above parameterisation.
We should emphasise that this parameterisation is such that if $y$ solves 
\PVI\ (with parameters $\theta$)
and $x$ is defined by \eqref{eq: p eq}
then the family of connections \eqref{eq: 3x3 conn} is isomonodromic as 
$t$ varies.
Indeed one may obtain a solution to Schlesinger's equations by also 
doing two quadratures as follows. 
(This amounts to varying the systems appropriately under the adjoint action of
the diagonal torus, which clearly only conjugates the monodromy.)
By construction the above parameterisation is transverse to the torus orbits. 
We will parameterise the torus orbits by
replacing $B_i$ above by $hB_ih^{-1}$ where $h=\diag(l,k,1)$ for parameters
$l,k\in\IC^*$.
One then finds the new residues $B_i$ solve Schlesinger's equations 
provided also
\beq\label{eq:dlogk}
\frac{d}{dt}\log k=\frac{\theta_4-1}{t(t-1)}(y-t)
\eeq
(as in \cite{JM81} p.445) and
\beq\label{eq:dlogl}
\frac{d}{dt}\log l=
\frac{\delta-1}{t(t-1)}\left(y-t-\frac{\de-\theta_4}{p}\right)
\eeq
where $p=x+\sum_1^3\theta_i/(y-a_i)$.
As a consistency check one can observe that the equations \eqref{eq:dlogk} and
\eqref{eq:dlogl} are exchanged by the transformation swapping $\mu_1$ and
$\mu_2$. Indeed the corresponding Okamoto transformation 
$(R_1R_2R_3)R_5(R_1R_2R_3)$ maps $y$ to $y-\frac{\de-\theta_4}{p}$ and changes
$\theta_4$ into $\de$.
\end{rmk}

\begin{rmk}   \label{rmk: simpler Bis}
The parameterisation of the $3\times 3$ Fuchsian systems given in the previous
remark is tailored so that one can see how the Okamoto transformation $R_5$
arises and see the relation to Schlesinger's equations (i.e. one may do 
the two quadratures to obtain a Schlesinger solution). 
However, when written out explicitly, the matrix entries are complicated
rational functions of $x,y,t,\lambda_i,\mu_i$ (the $2\times 2$ 
case in \cite{JM81} is already quite complicated).
If one is simply interested in writing down an isomonodromic family of
Fuchsian system (starting from a \PVI\ solution $y$) then one may conjugate
the above family of Fuchsian systems into a simpler form, as follows.
First, if we write each $B_i$ of the previous remark with respect to the basis
$\{f_i\}$, then $B_i$ will only have non-zero matrix entries in the $i$th row.
Then one can further conjugate 
by the diagonal torus to
obtain the following, simpler, explicit matrices:
\beq   \label{eq: Bi simple version}
B_1=\left(\begin{matrix}
\lambda_1& b_{12} & b_{13} \\
0&0&0\\
0&0&0
\end{matrix}\right),
\qquad
B_2=\left(\begin{matrix}
0&0&0 \\
b_{21}&\lambda_2&b_{23} \\
0&0&0
\end{matrix}\right),\qquad
B_3=\left(\begin{matrix}
0&0&0 \\
0&0&0 \\
b_{31}&b_{32}&\lambda_3
\end{matrix}\right)
\eeq
where
\begin{align*}
b_{12}&=
{ \lambda_1}\,\, - \mu_3y +(\mu_1 - xy)(y-1),
\ &
b_{32}&=(\mu_2-\lambda_2-b_{12})/t,
\\
b_{13}&=
  { \lambda_1}t-\mu_3y +(\mu_1 - xy)(y-t),
\ &
b_{23}&=(\mu_2-\lambda_3)t-b_{13},
\\
b_{21}&=
\lambda_2+\frac{\mu_3(y-t)-\mu_1(y-1)+x(y-t)(y-1)}{t-1},
\ &
b_{31}&=(\mu_2-\lambda_1-b_{21})/t.\\
\end{align*}
Thus if $y(t)$ solves \PVI\ (with parameters $\theta$ as in \eqref{thla})
and we define $x(t)$ via
\eqref{eq: p eq} and construct the matrices $B_i$ from the above formulae,
then the family of Fuchsian systems 
\beq   \label{eq: 3x3 system}
\frac{d}{dz}-\left(\frac{B_1}{z}+\frac{B_2}{z-t}+\frac{B_3}{z-1}\right)
\eeq
will be isomonodromic as $t$ varies,
since it is conjugate to a Schlesinger solution.
This seems to be the simplest way to write down explicit isomonodromic
families of rank three Fuchsian systems from \PVI\ solutions (an example will
be given in the following section). 
\end{rmk}

\end{section}

\begin{section}{Special solutions}

Another application of the $3\times 3$ Fuchsian representation of \PVI\ 
is that it
allows us to see new finite-branching solutions to \PVI.
The basic idea is that, due to \eqref{eq: f2 action},
if a Fuchsian system has finite linear monodromy group
then the solution to the isomonodromy equations, controlling its deformations,
will only have a finite number of branches.
For example this idea was used in the $2\times 2$ context
by Hitchin \cite{Hit-Poncelet,Hit-Octa} 
to find some explicit solutions with
dihedral, tetrahedral and octahedral linear monodromy groups.
(Also there are $5$ solutions in \cite{Dub95long,DubMaz00,Kitaev-dessins}
equivalent to solutions with icosahedral linear monodromy groups.)

One can also try to use the same idea in the $3\times 3$ context.
The first question to ask is: 
what are the possible finite monodromy groups of rank $3$
connections of the form \eqref{eq: 3x3 conn}?
Well (at least if $\lambda_i\not\in \IZ$), 
the local monodromies around $0,t,1$ will be conjugate to the exponentials of
the residues, which will be matrices of the form 
``identity $+$  rank one matrix'', i.e. they will be pseudo-reflections.
Moreover the finite groups generated by such pseudo-reflections, 
often called complex reflection groups, have been
classified by Shephard and Todd \cite{Shep-Todd}.
Looking at their list we immediately see that we get a richer class of finite
groups than the finite subgroups of $\SL_2(\IC)$, and so expect to get new
\PVI\  solutions.

For example the smallest non-real exceptional complex reflection group is the
Klein reflection group of order $336$ (which is a two-fold cover of Klein's
simple group of
holomorphic automorphisms of Klein's quartic curve).
This leads to:

\begin{thm}[\cite{k2p}]   \label{thm: Klein}
The rational functions
\begin{equation}\notag
y=-{\frac { \left( 5\,{s}^{2}-8\,s+5 \right)  \left( 7\,{s}^{2}-7\,s+4\right) }
{ s \left( s-2 \right)  \left( s+1 \right) \left( 2\,s-1 \right) 
\left( 4\,{s}^{2}-7\,s+7 \right)  }},\qquad
t={\frac { \left( 7\,{s}^{2}-7\,s+4 \right) ^{2}}{{s}^{3} \left( 4\,{s}^
{2}-7\,s+7 \right) ^{2}}},
\end{equation}
constitute a genus zero solution to \PVI\  with $7$ branches and parameters
$\theta=(2,2,2,4)/7$.
It governs isomonodromic deformations of a rank $3$ Fuchsian connection
of the form \eqref{eq: 3x3 conn} with 
linear monodromy group isomorphic to the Klein reflection group
and parameters 
$\lambda_i=1/2$, 
$(\mu_1,\mu_2,\mu_3)=(3,5,13)/14.$ 
Moreover this solution 
is not equivalent to (or a simple deformation of) 
any solution with finite $2\times 2$
linear monodromy group.
\end{thm}

As an example application 
of the formulae of remark \ref{rmk: simpler Bis} it is now easy
to write
down the corresponding isomonodromic family of rank three Fuchsian systems
having monodromy equal to the Klein complex reflection group (we have
conjugated the resulting system slightly to make it easier to write).
The result is that for any $s$ such that $t(s)\ne 0,1,\infty$ the system
\eqref{eq: 3x3 system}, with $t=t(s)$ as in Theorem \ref{thm: Klein},
has monodromy equal to the Klein reflection group, generated by reflections,
where the residues $B_i$ are given by \eqref{eq: Bi simple version}
with each $\lambda_i=1/2$ and 
\begin{align*}
b_{12}&=
{\frac {14\,{s}^{3}-21\,{s}^{2}+24\,s-22}
{21 s \left( 4\,{s}^{2}-7\,s+7 \right) }},&
b_{13}&={\frac {22\,{s}^{3}-24\,{s}^{2}+21\,s-14}
{21(7\,{s}^{2}-7\,s+4)}},\\
b_{21}&=
{\frac {14\,{s}^{3}-21\,{s}^{2}+24\,s+5}
{ 21\left( s-1 \right)  \left( 4\,{s}^{2}-s+4\right) }},&
b_{23}&=
{\frac {22\,{s}^{3}-42\,{s}^{2}+39\,s-5}{21(7\,{s}^{2}-7\,s+4)}},\\
b_{31}&=
{\frac {14-21\,s+24\,{s}^{2}+5\,{s}^{3}}{ 21\left( s-1 \right)
 \left( 4\,{s}^{2}-s+4 \right) }},&
b_{32}&=
{\frac {{22-42\,{s}+39\,s^2-5\,s^3}}
{21s \left( 4\,{s}^{2}-7\,s+7 \right) }}.
\end{align*}
Observe that 
$$t=
\frac { \left( 7\,{s}^{2}-7\,s+4 \right) ^{2}}
{{s}^{3} \left( 4\,{s}^{2}-7\,s+7 \right) ^{2}}=
1-{\frac { \left( 4\,{s}^{2}-s+4 \right) ^{2} \left( s-1 \right) ^{3}}{{s}^{3}
 \left( 4\,{s}^{2}-7\,s+7 \right) ^{2}}}$$
so that the matrix entries of the the residues 
$B_i$ are all nonsingular whenever
$t(s)\ne 0,1,\infty$.
(Up to conjugation, at the value $s=5/4$ this system
equals that of \cite{k2p} Corollary 31 although
there is a typo just before (p.200 \cite{k2p}) in that 
the values of $b_{23}b_{32}=\tr(B_2B_3)$ and 
$b_{13}b_{31}=\tr(B_1B_3)$ have been swapped.)

Unfortunately most of the other three-dimensional complex reflection groups 
do not seem to lead to new solutions of \PVI. However the largest 
exceptional complex reflection group does give new solutions. In this case the
group is the Valentiner reflection group of order $2160$ (which is a $6$-fold
cover of the group $A_6$ of even permutations of six letters).
Now one finds there are three inequivalent solutions that arise, all
of genus one.
(Choosing the linear monodromy representation amounts to choosing a triple of
generating reflections, and in this case there are three inequivalent triples
that can be chosen.)

\begin{thm}[\cite{icosa}]
There are three inequivalent triples of reflections generating the Valentiner
complex reflection group.
The \PVI\ solutions governing the isomonodromic deformations of the
corresponding Fuchsian systems are all of genus one. 
They have
$15, 15, 24$ branches and parameters
$$(\mu_1,\mu_2,\mu_3)=(5,11,29)/30,\qquad(5,17,23)/30,\qquad(2,5,11)/12,$$
respectively (with all $\lambda_i=1/2$). 
The explicit solutions appear in \cite{icosa}.
\end{thm}

Somewhat surprisingly when pushed down to the equivalent $2\times 2$
perspective these solutions all correspond to Fuchsian systems with
linear monodromy generating the binary icosahedral group, 
and they are not equivalent to any of the $5$
solutions already mentioned.
(The $3$ icosahedral solutions of Dubrovin and Mazzocco 
\cite{Dub95long,DubMaz00}, with $10, 10, 18$ branches respectively
do fit into this framework and correspond to the 
three inequivalent choices of
generating reflections of the icosahedral reflection
group, cf. also \cite{k2p} pp.181-183.) 

This led to the question of seeing what other such `icosahedral solutions' 
might
occur (e.g. is the $24$ branch solution the largest?).
The classification was carried out in \cite{icosa}.
(Another motivation was to find other interesting  
examples on which to apply the machinery
of \cite{Jimbo82, k2p} to construct explicit solutions.) 
At first glance one finds there is a huge number of such linear
representations; one is basically counting the number of conjugacy classes 
of triples of generators
of the binary icosahedral group, and an old formula of Hall \cite{Hall36}
says there are $26688$.
However this is drastically reduced if we agree to identify solutions if they
are related by Okamoto's affine $F_4$ action (since after all there is a
simple algebraic procedure to relate any two equivalent solutions, using the
formulae for the Okamoto transformations).

\begin{thm}[see \cite{icosa}]
There are exactly $52$ equivalence classes of solutions to \PVI\ having linear
monodromy group equal to the binary icosahedral group.

$\bullet$ 
The possible genera are: $0,1,2,3,7$, and the largest solution has $72$
branches.

$\bullet$ The first $10$ classes correspond to the ten icosahedral entries on
Schwarz's list of algebraic solutions to the hypergeometric equation,

$\bullet$ The next $9$ solutions have less than $5$ branches and are simple
deformations of known (dihedral, tetrahedral or octahedral) solutions,

\noindent
The remaining $33$ solutions are all now known explicitly, namely there are: 

$\bullet$ The $5$ already mentioned of Dubrovin, Mazzocco and
Kitaev in \cite{Dub95long,DubMaz00,Kitaev-dessins}, 

$\bullet$ The $20$ in \cite{icosa} including the three Valentiner
solutions, and 

$\bullet$ The $8$ in \cite{ipc}, constructed 
out of previous solutions  via quadratic transformations.
\end{thm}

In particular all of the icosahedral solutions with more than $24$ branches
(and in particular all the icosahedral solutions with genus greater than
one) were obtained from earlier solutions using quadratic 
transformations, so in this sense the $24$ branch Valentiner 
solution is the largest `independent' icosahedral solution (it
was certainly the hardest to construct).

The main idea in the classification was to sandwich the equivalence classes 
between two other, more easily computed, equivalence relations (geometric and
parametric equivalence), which in this case turned out to coincide.
A key step was to understand the relation between the linear monodromy data of
Okamoto-equivalent solutions, for which the geometric description in
Theorem \ref{thm: R5 from schles} of the 
transformation $R_5$ was very useful (see also \cite{IIS}). 

Examining the list of icosahedral solutions carefully it turns out that there
is one solution which is ``generic'' in the sense that its parameters lie on
none of the reflection hyperplanes of the $F_4$ or $D_4$ affine Weyl groups. 
This is closely related to the fact that the icosahedral rotation group $A_5$
has four non-trivial conjugacy classes: one can choose a triple of pairwise
non-conjugate elements generating $A_5$ whose product is in the fourth
non-trivial class. Viewing this triple as a representation of the fundamental
group of a four-punctured sphere and choosing a lift to $\SL_2(\IC)$
arbitrarily, gives the monodromy data of a Fuchsian system with such
generic parameters.

\begin{plaincor}[\cite{icosa}]
There is an explicit algebraic solution to the sixth Painlev\'e equation whose
parameters lie on none of the reflecting hyperplanes of Okamoto's affine $F_4$
(or $D_4$) action.
\end{plaincor}

This contrasts for example 
with the Riccati solutions whose parameters always lie on an
affine $D^-_4$ hyperplane (and needless to say no other explicit 
generic solutions
are currently known).

One can also carry out the analogous classification for the tetrahedral and
octahedral groups, and this led to five new octahedral solutions.
In more detail:

\begin{thm}[see \cite{octa}]
There are exactly $6$ (resp. $13$) 
equivalence classes of solutions to \PVI\ having linear
monodromy group equal to the binary tetrahedral (resp. octahedral) group.

$\bullet$ The first two solutions of each type correspond to the two
entries of the same type on
Schwarz's list of algebraic solutions to the hypergeometric equation,

$\bullet$ The next solutions (with less than $5$ branches) were
previously found by Hitchin \cite{Hit-Poncelet,Hit-Octa} and
Dubrovin \cite{Dub95long} (up to equivalence/simple deformation),

$\bullet$
A six-branch genus zero tetrahedral solution and two genus zero octahedral
solutions (with $6$ and $8$ branches resp.)
were found by Andreev and Kitaev \cite{And-Kit-CMP, Kitaev-dessins}, 

$\bullet$ All the solutions have genus zero except for one $12$ branch 
octahedral solution of genus one. The largest octahedral solution has $16$
branches and is currently the largest known genus zero solution.
\end{thm}

\end{section}

\renewcommand{\baselinestretch}{1}              %
\normalsize
\bibliographystyle{amsplain}    \label{biby}
\bibliography{../thesis/syr}    

\providecommand{\bysame}{\leavevmode\hbox to3em{\hrulefill}\thinspace}
\providecommand{\MR}{\relax\ifhmode\unskip\space\fi MR }
\providecommand{\MRhref}[2]{%
  \href{http://www.ams.org/mathscinet-getitem?mr=#1}{#2}
}
\providecommand{\href}[2]{#2}
\begin{thebibliography}{10}

\bibitem{And-Kit-CMP}
F.~V. Andreev and A.~V. Kitaev, \emph{Transformations {$RS\sp 2\sb 4(3)$} of
  the ranks {$\leq 4$} and algebraic solutions of the sixth {P}ainlev\'e
  equation}, Comm. Math. Phys. \textbf{228} (2002), no.~1, 151--176.

\bibitem{BJL79}
W.~Balser, W.B. Jurkat, and D.A. Lutz, \emph{{B}irkhoff invariants and
  {S}tokes' multipliers for meromorphic linear differential equations}, J.
  Math. Anal. Appl. \textbf{71} (1979), 48--94.

\bibitem{nabh}
O.~Biquard and P.~P. Boalch, \emph{Wild non-abelian {H}odge theory on curves},
  Compositio Math. \textbf{140} (2004), no.~1, 179--204.

\bibitem{Birman}
J.~S. Birman, \emph{Braids, links, and mapping class groups}, Princeton Univ.
  Press, Princeton, N.J., 1974.

\bibitem{icosa}
P.~P. Boalch, \emph{The fifty-two icosahedral solutions to {P}ainlev\'e {VI}},
  J. Reine Angew. Math., to appear, (math.AG/0406281, v.7).

\bibitem{ipc}
\bysame, \emph{{Higher genus icosahedral {P}ainlev\'e curves}},
  math.AG/0506407.

\bibitem{saqh}
\bysame, \emph{Quasi-{H}amiltonian geometry of meromorphic connections},
  math.DG/0203161.

\bibitem{octa}
\bysame, \emph{Some explicit solutions to the {R}iemann--{H}ilbert problem},
  math.DG/0501464.

\bibitem{smid}
\bysame, \emph{{S}ymplectic manifolds and isomonodromic deformations}, Adv. in
  Math. \textbf{163} (2001), 137--205.

\bibitem{bafi}
\bysame, \emph{G-bundles, isomonodromy and quantum {W}eyl groups}, Int. Math.
  Res. Not. (2002), no.~22, 1129--1166, math.DG/0108152.

\bibitem{pecr}
\bysame, \emph{Painlev\'e equations and complex reflections}, Ann. Inst.
  Fourier \textbf{53} (2003), no.~4, 1009--1022.

\bibitem{k2p}
\bysame, \emph{From {K}lein to {P}ainlev\'e via {F}ourier, {L}aplace and
  {J}imbo}, Proc. London Math. Soc. \textbf{90} (2005), no.~3, 167--208.

\bibitem{BbkLie}
N.~Bourbaki, \emph{{G}roupes et alg\`ebres de {L}ie. {C}hapitres 4,5 et 6},
  Masson, Paris, 1981.

\bibitem{conway-smith}
J.~H. Conway and D.~A. Smith, \emph{On quaternions and octonions: their
  geometry, arithmetic, and symmetry}, A K Peters Ltd., Natick, MA, 2003.
  \MR{2004a:17002}

\bibitem{Dett-Reit-p6}
M.~Dettweiler and S.~Reiter, \emph{Painlev\'e equations and the middle
  convolution}, preprint, October 1 2004.

\bibitem{Dub95long}
B.~Dubrovin, \emph{Geometry of 2{D} topological field theories}, Integrable
  Systems and Quantum Groups (M.Francaviglia and S.Greco, eds.), vol. 1620,
  Springer Lect. Notes Math., 1995, pp.~120--348.

\bibitem{DubMaz00}
B.~Dubrovin and M.~Mazzocco, \emph{Monodromy of certain {P}ainlev\'e-{V}{I}
  transcendents and reflection groups}, Invent. Math. \textbf{141} (2000),
  no.~1, 55--147. \MR{2001j:34114}

\bibitem{Hall36}
P.~Hall, \emph{The {E}ulerian functions of a group}, Quart. J. Math. Oxford
  Ser. 7 (1936), 134--151.

\bibitem{Hit-Poncelet}
N.~J. Hitchin, \emph{Poncelet polygons and the {P}ainlev\'e equations},
  Geometry and analysis (Bombay, 1992), Tata Inst. Fund. Res., Bombay, 1995,
  pp.~151--185. \MR{97d:32042}

\bibitem{Hit-Octa}
\bysame, \emph{A lecture on the octahedron}, Bull. London Math. Soc.
  \textbf{35} (2003), 577--600.

\bibitem{IIS}
M.~Inaba, K.~Iwasaki, and M.-H. Saito, \emph{B\"acklund transformations of the
  sixth {P}ainlev\'e equation in terms of {R}iemann-{H}ilbert correspondence},
  I. M. R. N. (2004), no.~1, 1--30, math.AG/0309341.

\bibitem{Jimbo82}
M.~Jimbo, \emph{Monodromy problem and the boundary condition for some
  {P}ainlev\'e equations}, Publ. Res. Inst. Math. Sci. \textbf{18} (1982),
  no.~3, 1137--1161. \MR{85c:58050}

\bibitem{JM81}
M.~Jimbo and T.~Miwa, \emph{Monodromy preserving deformations of linear
  differential equations with rational coefficients {II}}, Physica 2D (1981),
  407--448.

\bibitem{JMU81}
M.~Jimbo, T.~Miwa, and Kimio Ueno, \emph{Monodromy preserving deformations of
  linear differential equations with rational coefficients {I}}, Physica 2D
  (1981), 306--352.

\bibitem{Katz-rls}
N.~M. Katz, \emph{Rigid local systems}, Annals of Mathematics Studies, vol.
  139, Princeton University Press, Princeton, NJ, 1996. \MR{MR1366651
  (97e:14027)}

\bibitem{Kitaev-dessins}
A.~V. Kitaev, \emph{Dessins d'enfants, their deformations and algebraic the
  sixth {P}ainlev\'e and {G}auss hypergeometric functions}, nlin.SI/0309078,
  v.3.

\bibitem{Mal-imd1long}
B.~Malgrange, \emph{Sur les deformations isomonodromiques. {I}. singularites
  regulieres}, S\'eminaire E.N.S. Math\'ematique et Physique (Boston)
  (L.~Boutet~de Monvel, A.~Douady, and J.-L. Verdier, eds.), Progress in Math.,
  vol.~37, Birkh\"auser, 1983, pp.~401--426.

\bibitem{NY-so8}
M.~Noumi and Y.~Yamada, \emph{A new {L}ax pair for the sixth {P}ainlev\'e
  equation associated with $\hat{\mathfrak{so}}(8)$}, Microlocal Analysis and
  Complex {F}ourier Analysis (K.~Fujita T.~Kawai, ed.), World Scientific, 2002.

\bibitem{OkaPVI}
K.~Okamoto, \emph{Studies on the {P}ainlev\'e equations. {I}. {S}ixth
  {P}ainlev\'e equation {$P\sb {{\rm VI}}$}}, Ann. Mat. Pura Appl. (4)
  \textbf{146} (1987), 337--381. \MR{88m:58062}

\bibitem{Schles}
L.~Schlesinger, \emph{\"{U}ber eine {K}lasse von {D}ifferentialsystemen
  beliebiger {O}rdnung mit feten kritischen {P}unkten}, J. f\"ur Math.
  \textbf{141} (1912), 96--145.

\bibitem{Shep-Todd}
G.~C. Shephard and J.~A. Todd, \emph{Finite unitary reflection groups},
  Canadian J. Math. \textbf{6} (1954), 274--304.

\bibitem{Simpson-nabh}
C.~T. Simpson, \emph{Nonabelian {H}odge theory}, Proceedings of the
  International Congress of Mathematicians, Vol.\ I, II (Kyoto, 1990) (Tokyo),
  Math. Soc. Japan, 1991, pp.~747--756.

\bibitem{Sim94ab}
C.T. Simpson, \emph{Moduli of representations of the fundamental group of a
  smooth projective variety, {I, II}}, Publ. Math. I.H.E.S. \textbf{79, 80}
  (1994), 47--129, 5--79.

\bibitem{watanabePVI}
H.~Watanabe, \emph{Birational canonical transformations and classical solutions
  of the sixth {P}ainlev\'e equation}, Ann. Scuola Norm. Sup. Pisa Cl. Sci. (4)
  \textbf{27} (1998), no.~3-4, 379--425.

\end{thebibliography}
\end{document}